\documentclass{emsprocart}
\usepackage[T1]{fontenc}
\usepackage{lmodern}
\usepackage{a4wide}
\usepackage{verbatim}
\usepackage{tikz}
\usepackage[english]{babel}

\def\balpha{{a}}

\newcounter{noalgo}[section]
\newdimen\indentalgo
\newdimen\indentalgodec\indentalgo=0.0mm\indentalgodec=10mm

\newcommand{\If}{\advance\indentalgo by \indentalgodec {\bf if }}

\newcommand{\For}{\global\advance\indentalgo by \indentalgodec {\bf for }}

\newcommand{\Endindent}{\global\advance\indentalgo by -\indentalgodec}

\newdimen\decalage \decalage=0.5cm
\newcounter{algo} 

\newcommand{\PP}{\mathbb P}

\let\set\mathbb
\def\<<{\leavevmode
  \raise0.28ex\hbox{$\scriptscriptstyle\langle\!\langle$}\nobreak
  \hskip -.6pt plus.3pt minus.2pt\,}
\def\>>{\,\nobreak\hskip -.6pt plus.3pt minus.2pt
  \raise0.28ex\hbox{$\scriptscriptstyle\rangle\!\rangle$}}

\def\Aut{\mathop{\rm{Aut}}\nolimits }
\def\Gal{\mathop{\rm{Gal}}\nolimits }
\def\Pic{\mathop{\rm{Pic}}\nolimits }
\def\GL{\mathop{\rm{GL}}\nolimits }
\def\SL{\mathop{\rm{SL}}\nolimits }
\def\ol{\overline}

\def\CC{{\set C}}
\def\HH{{\set H}}
\def\PU{{{\set P}^1}}
\def\QQ{{\set Q}}
\def\bQ{{\overline {\set Q}}}
\def\RR{{\set R}}
\def\ZZ{{\set Z}}
\def\cL{{\cal L}}

\newtheorem{subproposition}[subsubsection]{Proposition}
\newtheorem{theorem}[subsection]{Theorem}

\makeatletter
\newenvironment{eqn}{\refstepcounter{subsection}
$$}{\leqno{\rm(\thesubsection)}$$\global\@ignoretrue}
\makeatother  

\makeatletter 
\newenvironment{subeqn}{\refstepcounter{subsubsection}
$$}{\leqno{\rm(\thesubsubsection)}$$\global\@ignoretrue}
\makeatother 

\providecommand{\myproofname}{Proof}

\usepackage{enumitem}

\contact[Jean-Marc.Couveignes@math.u-bordeaux1.fr]{Jean-Marc
  Couveignes, 
Univ. Bordeaux, IMB, UMR 5251, F-33400 Talence, France.
CNRS, IMB, UMR 5251, F-33400 Talence, France.
INRIA, F-33400 Talence, France}
\contact[edix@math.leidenuniv.nl]{Mathematisch Instituut, 
Universiteit Leiden, Niels Bohrweg~1, 2333 CA Leiden, Nederland}

\title[Approximate computations with modular curves]
      {Approximate computations with modular curves}

\author[Jean-Marc Couveignes and Bas Edixhoven]{Jean-Marc
  Couveignes\thanks{Research  supported by  ANR (project ALGOL ANR-07-BLAN-0248) and by
    DGA ma{\^\i}trise de l'information.}
 \ and Bas Edixhoven}

\begin{document}
\selectlanguage{english}

\begin{abstract}
This article gives an introduction for mathematicians interested in
numerical computations in algebraic geometry and number theory to some
recent progress in algorithmic number theory, emphasising the key role
of approximate computations with modular curves and their
Jacobians. These approximations are done in polynomial time in the
dimension and the required number of significant digits. We explain
the main ideas of how the approximations are done, illustrating them
with examples, and we sketch some applications in number theory.
\end{abstract}

\selectlanguage{english}



\maketitle

\section{Introduction}\label{sec:intro}
The purpose of this article is to give an introduction to the main
results of the book \cite{book} and their  generalization
in the  PhD thesis \cite{Bruin1} and in  \cite{Bruin2}, as well as
some applications, and most of all to
explain the essential role played by \emph{approximate
  computations}. The intended reader is a mathematician interested in
\emph{numerical computations} in algebraic geometry or number theory.

The results concern fast algorithms in number theory and more
precisely, fast computation of Fourier coefficients of modular
forms. These coefficients, with Ramanujan's $\tau$-function as a
typical example, have deep arithmetic significance and are important
in various areas of mathematics, from number theory and algebraic
geometry to combinatorics and lattices.

The fastest previously known algorithms for computing these Fourier
coefficients took exponential time, except in some special cases. The
case of elliptic curves (Schoof's algorithm) was at the birth of
elliptic curve cryptography around 1985. The results mentioned above
give an algorithm for computing coefficients of modular forms in
polynomial time. For example, Ramanujan's $\tau(p)$ with $p$ a prime
number can be computed in time bounded by a fixed power of~$\log p$.

Such fast computation of Fourier coefficients is itself based on the
main result of the book: the computation, in polynomial time, of
Galois representations over finite fields attached to modular forms by
the Langlands program.

The computation of the Galois representations uses their realisation,
following Shimura and Deli\-gne, in the torsion subgroup of Jacobian
varieties of modular curves. The main challenge is then to perform the
necessary computations in time polynomial in the dimension of these
nonlinear algebraic varieties. Exact computations involving systems of
polynomial equations in many variables take exponential time. This is
avoided by numerical approximations with a precision that suffices to
derive exact results from them. Bounds for the required precision -- in
other words, bounds for the height of the rational numbers that
describe the Galois representation to be computed -- are obtained from
Arakelov theory.

This article is organised as follows. Sections~\ref{sec:afcaj}
and~\ref{sec:psi} are concerned with numerical methods used in the
context of complex algebraic curves and their Jacobian
varieties. Sections~\ref{sec:ctd1} and~\ref{sec:ctd2} describe how to
get exact results about torsion points on modular curves using these
numerical methods. Section~\ref{sec:ctd1} focuses on the genus $1$
curve $X_{11}$ while Section~\ref{sec:ctd2} deals with the general
modular curve $X_\ell$. As an application, Section~\ref{sec:avoq}
gives two examples of fast computation of coefficients of modular
forms: Ramanujan's $\tau$-function, and the classical sums of squares
problem.

\section{Algorithms for curves and  Jacobians}\label{sec:afcaj}

Let $X$ be a connected, smooth, projective algebraic curve over
the field $\CC$ of complex numbers. 
The set  $X(\CC)$ of complex
points of $X$ is a Riemann surface. Let $g$ be the genus of $X$ and
let $(\omega_k)_{1\le k\le g}$ be a basis for  the space of holomorphic
differentials on $X$. We fix a point   $b\in X(\CC)$ and we denote
by $Y_b$ the set of homotopy classes of paths on $X(\CC)$ starting at $b$.
The {\it universal cover}
$f_b : Y_b \rightarrow X(\CC)$  maps  every path to its end point.
The {\it fundamental group} $\pi_1(X(\CC),b)\subset Y_b$ is
 the subset of (homotopy classes of)
closed paths. It acts on~$Y_b$, with quotient~$X(\CC)$.
We have an integration map $\phi_b : Y_b\rightarrow \CC^g$ defined
by \[\phi_b(\gamma)=(\int_\gamma\omega_1, \ldots, \int_\gamma
\omega_g).\] The image of $\pi_1(X(\CC),b)$ by $\phi_b$ is a lattice
$\Lambda$ in $\CC^g$.  It is called the {\it lattice of periods}. It
is a free $\ZZ$-module of rank~$2g$.  The quotient $\CC^g/\Lambda$ is
a complex torus. It is  the set of complex points $J(\CC)$ on the
{\it Jacobian variety} $J$ of $X$.  The integration map $\phi_b :
Y_b\rightarrow \CC^g$ induces a map between the quotients
$X(\CC)\rightarrow J(\CC)$. This map is a morphism of varieties
$X\rightarrow J$. We call this morphism $\phi_b$ also. For every
positive integer $k$ we denote $\phi_b^k : X^k\rightarrow J$ the
morphism that maps $(P_1, \ldots, P_k)$ onto
$\phi_b(P_1)+\dots+\phi_b(P_k)$.  Since the image in $J$ does not
depend on the ordering on the points $P_j$, we write $X^{(k)}$ for the
$k$-th symmetric power of~$X$. We note that $X^{(k)}$ is the quotient
of $X^k$ by the action of the symmetric group. It is a nonsingular
variety.  
We define the morphism $\phi_b^{(k)} : X^{(k)}\rightarrow J$ that
maps $\{P_1, \ldots, P_k\}$ onto $\phi_b(P_1)+\dots+\phi_b(P_k)$. 
For $k=g$ the map $\phi_b^{(g)}$ is
birational and surjective.  It is not an isomorphism unless $g\leq 1$.
Its fibers are projective linear spaces, mostly (but not all) points.
A degree $g$ effective divisor $P= P_1+\dots+P_g$ is said to be {\it
  non-special} if the map $\phi^{(g)}_b$ is a local diffeomorphism at
$P$. Otherwise we say that $P$ is {\it special}. This definition does
not depend on the chosen origin $b$. The set of special effective
degree $g$ divisors is the singular locus of $\phi_b^{(g)}$.  All
these maps $\phi^{(k)}_b$ are called Abel-Jacobi maps. In particular
\[
\phi_b^{(g)}(\{P_1, \ldots, P_g\})=
\sum_{1\le j\le g}(\int_b^{P_j}\omega_k)_{k} \, \bmod \Lambda,
\]
where we can  integrate $\int_b^{P_j}\omega_k$  along any path between
$b$ and $P_j$, provided we keep the  same path for all $k$.
We can apply the Abel-Jacobi map to any divisor on $X$. We set 
$\phi_b(\sum_j e_jP_j)=\sum_je_j\phi_b^{(1)}(P_j)$. We note that for degree
zero divisors, the image does not depend on the origin $b$.
A  divisor is said to be {\it principal} if it is the divisor
of a non-zero meromorphic function on $X$. Two divisors are said
to be {\it linearly equivalent} when their difference is principal.
Any principal divisor has degree  zero.
A degree zero divisor is principal if and only its image by $\phi_b$
is zero. So the set $J(\CC)=\CC^g/\Lambda$ of complex points on the Jacobian
is canonically identified with the group $\Pic^0(X)$
of 
linear equivalence classes of degree zero divisors  on $X$.

We now list important algorithmic problems related to the Abel-Jacobi
map. We illustrate them on the simple example of the projective
curve $X$ with
equation
\begin{subeqn}\label{eq:X11}
Y^2Z-YZ^2=X^3-X^2Z.
\end{subeqn} This curve has genus $1$. We write $x=X/Z$
and $y=Y/Z$. The  unique (up to a multiplicative constant) holomorphic
differential on $X$ is \[\omega= \frac{dx}{2y-1}=\frac{dy}{x(3x-2)}.\]
We choose the point $b=[0:1:0]$ as origin  for the integration map.
For every computational problem we shall consider, we will also 
explain what can be proven when $X$ is a modular curve $X_\ell$
and 
$\ell$ (therefore $g$) 
tends to infinity. The definition of the modular curve $X_\ell$ is given
in Section~\ref{sec:ctd2}. See also textbooks  \cite{Diamond-Shurman, Stein}
where $X_\ell$ is often denoted $X_1(\ell)$.

\subsection{Computing the lattice of periods}\label{sec:lattice}

We first need a basis for the singular homology group $H_1(X(\CC),\ZZ)$.
 If $X$
is the genus one curve given by equation~(\ref{eq:X11}),
such a  basis can be deduced from the study of
the degree two map $x : X\rightarrow \PU$
that sends $(x,y)$ onto $x$ and $[0:1:0]$ to $\infty$.
This map is ramified at $\infty$ and the three roots of $4x^3-4x^2+1$. We lift
 a simple loop around $\infty$ and one of
these  three roots. We then lift
 a simple loop around $\infty$ and another   root.
We thus obtain two elements in   $H_1(X(\CC),\ZZ)$ that form
a basis for it. 

Integrating a differential along a path is easy. We express
the differential in terms of  local
coordinates. We then reduce  to integrating   converging
power series. We integrate term by term. 
In case $X$ is the curve given in equation~(\ref{eq:X11}), we obtain a
basis $(\Omega_1, \Omega_2)$ for the lattice $\Lambda$ of periods where
\begin{eqnarray*}
\Omega_1&=& 6.346046521397767108443973084,\\ \Omega_2&=&-3.173023260698883554221986542 + 1.458816616938495229330889613i.\end{eqnarray*}

These calculation are made e.g. using the \cite{pari} system.
\begin{verbatim}
>a1=0;a2=-1;a3=-1;a4=0;a6=0;
>X=[a1,a2,a3,a4,a6];X=ellinit(X);
>X.omega
[6.346046521397767108443973084, 
-3.173023260698883554221986542 + 1.458816616938495229330889613*I]
\end{verbatim}
When dealing with general modular curves, an explicit basis for both
the singular homology and the de~Rham cohomology is provided by the
theory of Manin symbols \cite{manin, merel, cremona, frey,
  Stein}. Computing (good approximations of) periods is then achieved
in time polynomial in the genus and the required accuracy
\cite{Couveignes2}. The practical side is described
in~\cite[\S6.3]{Bosman1}. Textbooks \cite{Cohen}, \cite[Chapter
  3]{cremona} give even faster techniques for genus $1$ curves, but we
shall not need them.

\subsection{Computing with divisor classes}\label{sec:divcla}

A degree zero divisor class can  be represented by  
a point in the torus $\CC^g/\Lambda=J(\CC)$. It can also
be represented by  a  divisor 
of the form \begin{subeqn}\label{eq:div}
P_1+\dots+P_g-gb
\end{subeqn}
in this class.  This latter representation is not always
 unique. It is however unique for most classes because 
$\phi_b^{(g)}$ is birational. 
The addition problem in this context is the following: 
given two degree $g$ effective  divisors
$P=P_1+\dots+P_g$ and $Q=Q_1+\dots+Q_g$, one would like to compute
a degree $g$ effective divisor
$R=R_1+\dots+R_g$ such that
the divisor class of $R-gb$
is the  sum of the divisor classes of $P-gb$ and $Q-gb$. 
So we look for $g$ complex points
$R_1$, \ldots, $R_g$ such that $P_1+\dots+P_g+Q_1+\dots+Q_g-2gb$
is linearly equivalent to $R_1+\dots+R_g-gb$.  This is achieved using the
Brill-Noether
algorithm \cite{BN,Vol}. This algorithm uses a complete linear space $\cL$
of forms or functions. 
This space should have dimension $\ge 2g+1$. For example, assuming
$g\ge 4$,
we may take for $\cL$ the space of all holomorphic quadratic  differential 
forms.
We compute once for all a basis for this space. Then  the Brill-Noether
algorithm alternates  several steps of two different natures. Sometimes
we are given a form (function) and we want to compute its divisor.
Sometimes we are given an effective  divisor $D$
and we want to compute
a basis for  the subspace  $\cL(-D)$ consisting of forms (functions) 
vanishing  at this divisor.

The first problem (finding zeros of a given form) can be reduced,
using a convenient coordinate system, to the following problem: given
a power series $f(z)=\sum_{k\ge 0} f_kz^k$ with radius of convergence
$\ge 1$, find approximations of its zeros in the disk $D(0,1/2)$ with
center $0$ and radius $1/2$.  It is clear (see
\cite[\S5.4]{Couveignes1}) that, for the purpose of finding zeros, one
can replace $f(z)$ by its truncation $\sum_{0\le k\le K}f_kz^k$ at a
not too large order $K$. We then reduce to the classical problem of
computing zeros of polynomials. A survey of this problem is given in
\cite[\S5.3]{Couveignes1}.

The second problem (finding the  subspace of functions vanishing at given
points) boils down to finding the kernel of the  matrix  having  entries
the values of the  functions in the chosen basis of $\cL$
at the given points. 

The only difficulty
then is to control the conditioning of these two  problems. This is done
in two steps. 
We first 
 prove \cite[\S5.4]{Couveignes1} that the zeros of a holomorphic function
on a closed disk are well conditioned unless this function is small
everywhere on this disk. We then
prove \cite[\S12.7]{Couveignes2} that the form we consider
cannot be small everywhere on any of  the charts we consider, unless
it has very small coordinates in the chosen  basis of $\cL$.

The resulting algorithm for computing in the group of divisor classes
of modular curves is polynomial time in the genus and the required
direct accuracy \cite[Theorem 12.9.1]{Couveignes2}. By {\it direct
  accuracy} we mean that the error is measured in the target space of
the integration map, namely the torus $\CC^g/\Lambda$.  Saying that
the direct accuracy is bounded from above by $\epsilon$ means that the
returned divisor $R'=R_1'+\dots+R_g'$ is such that
\[\phi_b (R'-R) = \phi^{(g)}_b(R')-\phi^{(g)}_b(R)\] is bounded from above
by $\epsilon$ for the maxnorm in $\CC^g$.
This does not necessarily imply that the $R_j$ are close to
the $R_j'$. Indeed, in case
$R=R_1+\dots+R_g$ is special, there exists a non-trivial
 linear pencil of divisors $R'$ such that $\phi_b^{(g)}(R')=\phi_b^{(g)}(R)$.
Controlling the distance between $R$ and $R'$ will only be possible 
in some cases.

In the special case when $X$ is the curve given by equation~(\ref{eq:X11})
 the map $\phi^{(1)}_b : X\rightarrow J$ is an isomorphism 
because the genus is $1$.
Computing with divisor classes is then very simple and
the Brill-Noether
algorithm takes a simple  form. The space $\cL$ consists
of all degree $1$ homogeneous forms,  and a basis for it is made of the
three projective coordinates  $X$, $Y$ and $Z$.
 Given $P$ and $Q$, one 
considers the unique projective line $\Delta_1$ through $P$ and $Q$.
In case $P=Q$ we take $\Delta_1$ to be the tangent to $X$ at $P$.
The line $\Delta_1$ meets $X$ at three points: $P$, $Q$ and a third
point that we call $S$. We consider the unique projective line 
$\Delta_2$ through $S$ and the origin $b$.
The line $\Delta_2$ meets $X$ at three points: $b$, $S$ and a third
point that we call $R$. On can easily check that $P+Q$
is linearly equivalent to $b+R$ or equivalently $P-b+Q-b$
is equivalent  to $R-b$.
The coordinates of $R$ can be computed using very simple formulae 
\cite[Chapter III]{silverman}.
We illustrate this using  the \cite{pari} system. We call $P$ the 
point $[0:0:1]$. We first compute $Q$ such that $Q-b$ is linearly equivalent 
to $2(P-b)$. We write $Q-b\equiv 2(P-b)$ using the $\equiv$ symbol
for linear equivalence. We then compute $R$ such that  $R-b\equiv P-b+Q-b
\equiv 3(P-b)$. We then compute $S$ such that  $S-b\equiv Q-b+R-b\equiv
5(P-b)$.

\begin{verbatim}
>P=[0,0];
>Q=elladd(X,P,P)
[1, 1]
>R=elladd(X,P,Q)
[1, 0]
>S=elladd(X,Q,R)
[0]
\end{verbatim}
The answer for $S$ means that $S$ is just the origin $b=[0:1:0]$.
So the divisor $P-b$ has order $5$ in the Picard group 
$\Pic(X)$, the group
of divisors modulo linear equivalence.

\subsection{The direct Jacobi problem}\label{sec:jaco}

Given a  divisor  on $X$ we want to compute its image
by $\phi_b$ in the complex torus $J(\CC)=\CC^g/\Lambda$.
It suffices to explain what to do when the divisor
consists of a single point $P$. For every $1\le k\le g$
we then have to compute $\int_b^P\omega_k$. So we integrate
$\omega_k$ along any path from $b$ to $P$.   We split the chosen
path in several pieces according to the various charts in 
our atlas for the Riemann
surface $X(\CC)$.  On every chart, the differentials $\omega_k$
can be expressed in terms of the local coordinate. We then 
reduce to computing integrals of the form $\int_0^{\frac{1}{2}}f(z)dz$
where $f(z)$ is holomorphic on the unit disk. Such an  integral
can be computed  term by term. When $X$ is a modular curve, we
have a convenient system of charts and a basis for $\cL$
consisting of forms
having small coefficients in their expansions at every chart. 
There is long standing tradition with stating and proving bounds for these
coefficients. It culminates with 
the so-called Ramanujan conjecture.
This conjecture was proved by
Deligne as a consequence of~\cite{Deligne1} and his proof of the
analog of the Riemann hypothesis in the Weil conjectures
in~\cite{Deligne2}.
In case $X$ is the  elliptic curve given by equation~(\ref{eq:X11})
we  take for $P$ the point $[0:0:1]$ and find
that \[\phi_b^{(1)}(P)=\int_b^P\omega = 2.538418608559106843377589234\bmod \Lambda.\]
This integral is computed  using the \cite{pari} system.
\begin{verbatim}
> ellpointtoz(X,[0,0])
2.538418608559106843377589234
\end{verbatim}
We notice that \[\phi_b^{(1)}(P)=\frac{2\Omega_1}{5}\bmod \Lambda.\] So $5(P-b)$ is a principal divisor
as already observed at the end of section~\ref{sec:jaco}.

\subsection{The inverse Jacobi problem}\label{sec:jacinv}

At this point we have two different ways of representing a degree
zero class of equivalence of divisors. We can be given a divisor
in this class like  the one in equation~(\ref{eq:div}). Such a divisor
will be called 
a {\it reduced divisor}. We can also be given a vector  in $\CC^g$
modulo the lattice of periods $\Lambda$.
It is of course very easy to compute with such vectors. We also
have seen in section~\ref{sec:divcla} how to compute with 
reduced divisors. So both  representations  are  convenient for
computational purposes. We also have seen in section~\ref{sec:jaco}
how to pass from a reduced divisor to the corresponding point
in the torus $\CC^g/\Lambda$ applying the Abel-Jacobi map.
We now consider the inverse problem: given a point $\alpha\bmod \Lambda$
in the torus $\CC^g/\Lambda$, find  some $P=P_1+\cdots+P_g$
such that the   reduced divisor 
$P-gb$  is mapped onto $\alpha \bmod \Lambda$ by $\phi_b$.

\paragraph{Using an iterative method}

We can try an iterative method like  the secant's method. We illustrate
the secant's method in case $X$ is the curve
given by equation~(\ref{eq:X11}) and 
\begin{subeqn}\label{eq:alpha}
\alpha=(\Omega_1+\Omega_2)/11=0.2884566600635348685656 
+ 0.1326196924489541117573i.
\end{subeqn}
 Starting from  $P_0=(50-50i,-223.147+547.739i)$ and 
$P_1=(20-20i,-54.587+137.965i)$ we obtain an approximation
up to $10^{-26}$ after eighteen iterations. 
We use the \cite{pari} system and declare a function
for the secant method.
\begin{verbatim}
>secant(alpha,P0,P1,K)=
{
local(f0,f1,x0,x1,x2,P2,P3);
for(k=1,K,
f0=ellpointtoz(X,P0)-alpha;f1=ellpointtoz(X,P1)-alpha;
x0=P0[1];x1=P1[1];
x2=x1-f1*(x1-x0)/(f1-f0);
P2=[x2,ellordinate(X,x2)[1]];P3=[x2,ellordinate(X,x2)[2]];
if(abs(P2[2]-P0[2])> abs(P3[2]-P0[2]) ,P2=P3,);
P0=P1;P1=P2;
);
return(P2);
}
\end{verbatim}
The four parameters of this function are the target  point in $\CC/\Lambda$,
the two initial approximate 
values of $P$, and the number of iterations.
We then type
\begin{verbatim}
>alpha=(omega1+omega2)/11;
x0=50-50*I;x1=20-20*I;
P0=[x0,ellordinate(X,x0)[2]];P1=[x1,ellordinate(X,x1)[2]];
secant(alpha,P0,P1,18)
\end{verbatim}
Below are the results
of
iterations
$14$ to $18$. We only give the values taken by the $x$-coordinate.
\begin{verbatim}
6.796891402429021881380876803 - 7.525836023544396684018482041i
6.796539495414535904114103146 - 7.525907619429540863361002543i
6.796539142100022043003057330 - 7.525908029913269174706910680i
6.796539142094915910541452272 - 7.525908029899464322147329306i
6.796539142094915911068237206 - 7.525908029899464321854796862i
\end{verbatim}

\paragraph{The continuation method}
Iterative methods only work if the starting approximation
is close enough to the actual solution. 
Such an initial approximation can be provided  by the
solution of a different though close inverse problem.
Coming back to  our example, we will start from any
point on $X$. Say $P_0=(0,0)$. We compute the image
$\alpha_0\bmod \Lambda$ of $P_0$  by the integration map.
We then choose any $P_{-1}$ that is close enough to $P_0$.
\begin{verbatim}
>P0=[0,0];
alpha0=ellpointtoz(X,P0);
Pm1=[0.1,ellordinate(X,0.1)[2]];
\end{verbatim}
We now move slowly  from $\alpha_0$ to $\alpha$.
We set $\alpha_1=\alpha_0+0.1(\alpha-\alpha_0)$
and we solve the inverse problem for $\alpha_1$
using the secant's method with initial values
$P_{-1}$ and $P_0$.
\begin{verbatim}
>P1=secant(alpha0+0.1*(alpha-alpha0),Pm1, P0,5)
[0.218773824415936734050679268 - 0.0122309960881052801981765895*I,
0.0388323642082357612959944279 - 0.00390018046133107189481433241*I]
\end{verbatim}
We now set $\alpha_2=\alpha_0+0.2(\alpha-\alpha_0)$
and we solve the inverse problem for $\alpha_2$
using the secant's method with initial values
$P_0$ and $P_1$.
\begin{verbatim}
>P2=secant(alpha0+0.2*(alpha-alpha0),P0, P1,5)
[0.410237833586311839505201998 - 0.0205989424813431290064696558*I,
0.111775424533436210193603161 - 0.00838376796781394064004855129*I]
\end{verbatim}
We continue until we reach $\alpha$
\begin{verbatim}
>P3=secant(alpha0+0.3*(alpha-alpha0),P1, P2,5);
                    ...
P9=secant(alpha0+0.9*(alpha-alpha0),P7, P8,5);
P10=secant(alpha,P8, P9,10)
[6.796539142094915911068237205 - 7.525908029899464321854796861*I,
-8.056577776742775028742861296 + 30.05694612451787404370259256*I]
\end{verbatim}
This continuation method is very  likely to succeed provided
the integration map has a nice local behaviour all along the path
from $\alpha_0$ to $\alpha$. This is how practical computations have been
realised in \cite{Bosman1} for modular curves.
It is however  difficult to prove
that this method works because the integration map $\phi^{(g)}_b$ has
a singular locus as soon as $g>1$,  and we do not know how to provably
and efficiently find a path from $\alpha_0\bmod \Lambda$
to $\alpha\bmod \Lambda$ that keeps away from
the singular locus.

\section{Provably solving the inverse Jacobi problem}\label{sec:psi}

We have presented  in section~\ref{sec:jacinv} a heuristic algorithm
for the inverse Jacobi problem. This algorithm is based on  
continuation. It seems difficult to prove it however because that would
require a good control on the singular locus of the Jacobi map.
In this section we present the algorithm introduced in 
\cite{Couveignes2}. This algorithm only requires a good control
of the Jacobi map locally at a chosen divisor in $X^{(g)}$.
This is a much weaker condition and it is satisfied  for modular curves.
An important feature of this algorithm is the use of fast exponentiation
rather than  continuation. The principle of fast 
exponentiation is recalled in section~\ref{sec:fe}.   The algorithm
for the inverse Jacobi problem itself is given in section~\ref{sec:jil}.
Section~\ref{sec:mp} sketches the proof of this algorithm.
Proving in this context means proving the existence of 
a Turing  machine that returns a correct answer
in a given time. One has to prove both the correctness
of the result and a bound for the running time.
This bound here will be polynomial in the genus of the
curve and the required accuracy of the result.

\subsection{Fast exponentiation in groups}\label{sec:fe}

Assume we are given a  group $G$. The group law in $G$ will
be denoted multiplicatively.
 We assume that $G$ is  {\it computational}.
This  means that we know how to represent elements in $G$, how to compare
two given  elements, how to invert a given element,
and how to multiply two given elements.

The exponentiation problem in $G$ is the following: 
we are given an element $g$ in $G$ and
 an  integer $e\ge 2$, and we want to compute $g^e$ as an element
in $G$.
A first possibility  would be to set  $a_1=g$ and 
to compute $a_k=a_{k-1}\times g$ 
for $2\le k\le e$.
This requires $e-1$ multiplications in $G$.
It is well known, however, that we can do much better. We write
the  expansion of $e$ in base $2$,
\[e=\sum_{0\le k\le K}\epsilon_k2^k,\]
and we set $b_0=g$ and
$b_{k}=b_{k-1}^2$ for $1\le k\le K$. We then notice that
\[g^e=  \prod_{0\le k\le K}b_k^{\epsilon_k}.\]
So we can compute $g^e$ at the expense of a constant times
$\log e$ operations in $G$.
The algorithm above is called {\it fast exponentiation} and it
admits many variants and improvements \cite{gor}. Its first known occurrence
dates back to Pi\.{n}gala's Chandah-s{\^u}tra (before -200).
See \cite[I,13]{hist}.

\subsection{Solving the Jacobi inverse problem by linear algebra}\label{sec:jil}

Recall that   we have two different ways of representing 
an equivalence class of divisors of degree zero:
reduced divisors
or classes in the torus $\CC^g/\Lambda$.
We have seen that both  models are computational.
The Abel-Jacobi map $\phi^{(g)}_b : X^{(g)}\rightarrow \CC^g/\Lambda$
is computational also. We want to invert it (although we know
it is not quite injective). More precisely we assume we are given
some $\alpha$ in $\CC^g$ and we look for a degree $g$ effective
divisor on $X$ such that $\phi^{(g)}_b(P)=\phi_b (P-gb)=\alpha\bmod \Lambda$.
It seems difficult to prove the heuristic methods given
in section~\ref{sec:jacinv} for this purpose. So we present here a
 variant for which 
we can give a proof, at least when $X$ is a modular curve $X_\ell$.
We illustrate this method in the case where $X$ is the curve
given in equation~(\ref{eq:X11}). We still aim at the $\alpha$
given in equation~(\ref{eq:alpha}).

We need  a non-special effective divisor $P_0$  of degree
$g$. Since $g=1$ we can take any point  on $X$. For example
$P_0=(0,0)$. We note that the affine coordinate  $x$
is a local parameter at $P_0$.
We  choose a small real number $\epsilon$.
The smaller $\epsilon$ the better the precision of the final
result. Here we choose $\epsilon = 0.0001$. 
We consider two points $P_1$ and $P_2$ that are very 
close to $P_0$. The first point $P_1$ is obtained by
adding $\epsilon$ to the $x$-coordinate of $P_0$.
The second  point $P_2$ is obtained by
adding $\epsilon i$ to the $x$-coordinate of $P_0$.
\begin{verbatim}
P0=[0,0];
P1=[0.0001,ellordinate(X,0.0001)[2]];
P2=[0.0001*I,ellordinate(X,0.0001*I)[2]];
\end{verbatim}
We now compute the image $\alpha_1\bmod \Lambda$  of $P_1-P_0$ 
by the Abel-Jacobi
map. We also compute the image $\alpha_2\bmod \Lambda$ of  $P_2-P_0$.
We note that $\alpha_1\bmod \Lambda$ is very close to $0\in \CC/\Lambda$.
This is because $P_0$ and $P_1$ are close.
We assume that  $\alpha_1$ is the smallest
complex number
in its class modulo  $\Lambda$.
We make the same assumption for $\alpha_2$. 
Then $\alpha_1$ and $\alpha_2$ are two small complex 
numbers, and they form an $\RR$-basis of $\CC$.
This is  because the integration map $\phi^{(g)}_b$
is a local diffeomorphism at $P_0$ (or equivalently $P_0$ 
is a non-special divisor) and $\epsilon$ has been chosen
small enough.
\begin{verbatim}
alpha1=ellpointtoz(X,P1)-ellpointtoz(X,P0);
alpha2=ellpointtoz(X,P2)-ellpointtoz(X,P0)-omega1-omega2;
\end{verbatim}
Recall that our target in the torus $\CC/\Lambda$
is $\alpha \bmod \Lambda$ where $\alpha$ is the complex
number given in equation~(\ref{eq:alpha}). 
So we compute the two real coordinates of  $\alpha$ in the
basis $(\alpha_1,\alpha_2)$. 
\begin{verbatim}
>M=[real(alpha1), real(alpha2); imag(alpha1), imag(alpha2)]; 
coord=M^(-1)*[real(alpha),imag(alpha)]~
[-2884.566581407009845250155464, -1326.196933330853847302268151]~
\end{verbatim}
We deduce that $\alpha$ is very close to $\alpha' = -2884\alpha_1-1326\alpha_2$.
And the class $\alpha'\bmod \Lambda$ is the image by $\phi_b$  of
$-2884(P_1-P_0)-1326(P_2-P_0)$. The linear equivalence
class of the latter divisor is therefore a good
approximation for our problem. There remains to compute a {\it reduced
divisor} $P-gb$ in this class using the methods presented 
in section~\ref{sec:divcla}. Since the integers $2884$ and $1326$ are
rather big, we  use the fast exponentiation algorithm 
presented in section~\ref{sec:fe}.
\begin{verbatim}
>coord=truncate(coord)
[-2884, -1326]~
>D1=ellsub(X,P1,P0);D2=ellsub(X,P2,P0);
P=elladd(X,ellpow(X,D1,coord[1]),ellpow(X,D2,coord[2]))
[6.798693122986621316758396123 - 7.528977879167267357619566769*I, 
-8.059779911380488392224788509 + 30.07437308400090422713306570*I]
\end{verbatim}
We now check that the image of $P-P_0$ by $\phi_b$ is close to $\alpha$
\begin{verbatim}
>ellpointtoz(P)
0.2884000018811813146007079855 + 0.1325999988977252987328424662*I
>alpha
0.2884566600635348685656351402 + 0.1326196924489541117573536012*I
\end{verbatim}
For a better approximation we should start with  a smaller $\epsilon$.

\subsection{Matter of proof}\label{sec:mp}

The main concern when proving  the algorithm in section~\ref{sec:jil}
is to prove that we can find an initial divisor $P_0$ that is non-special.
In fact we must  guarantee a quantified version of this non-speciality
condition. The differential of $\phi^{(g)}_b$ at $P_0$ should be non singular
and its norm should not be too small. We can prove that such a condition
holds true for modular curves \cite[\S12.6.7]{Couveignes2} 
because we have a very sharp description
of these curves in the neighbourhood of the points called {\it cusps}.
As a consequence we prove \cite[Theorem 12.10.5]{Couveignes2} that the
inverse Jacobi problem for modular curves
can be solved in deterministic polynomial time 
in the genus and the required {\it direct accuracy}. Recall that
 {\it direct accuracy} means that the error is measured in the 
target space $\CC^g/\Lambda$.
The main difference between the algorithm in this section
and the one in section~\ref{sec:jacinv} is that we only need here to control
the local behaviour of $\phi^{(g)}_b$ at $P_0$ while the algorithm
in section~\ref{sec:jacinv} requires that the map $\phi^{(g)}_b$ be non-singular
above the whole path from $\alpha_0$  to $\alpha$.

In some cases it will be desirable  to control the {\it inverse error}
that is the error on the output divisor $P$ in $X^{(g)}$. This will
be possible when we can prove that $\phi^{(g)}: X^{(g)}\rightarrow
J$ is a local diffeomorphism at $P$ (that is $P$ is  non-special).
We will also need a lower bound for the norm of the differential
of $\phi^{(g)}$ at $P$. Such a lower bound  can be provided by arithmetic.

\section{Computing  torsion points I}\label{sec:ctd1}

In this and the next 
section we will assume that $X$ is a modular curve and $\ell$
a prime number. We will be interested in $\ell$-torsion points in 
the torus
$J(\CC)=\CC^g/\Lambda$. 
A point
\[\balpha = \alpha\bmod \Lambda\] is an $\ell$-torsion point
if and only if $\alpha$ lies in $\frac{1}{\ell}\Lambda$.
So the $\ell$-torsion subgroup
of $J(\CC)$ is 
$\frac{1}{\ell}\Lambda / \Lambda$ and it has  cardinality $\ell^{2g}$.
This group is also denoted $J[\ell]$.

 Some of these torsion points
 carry important arithmetic information. The values
taken by algebraic functions at these points generate interesting number
fields.  We want to  compute these fields. 
In this section
we will focus on a special case. We will assume that $X$ is the genus
$1$ curve given in equation~(\ref{eq:X11}) and $\ell=11$.
A  more general situation  will be studied  in the next
section~\ref{sec:ctd2}. We notice that the curve in equation~(\ref{eq:X11})
is indeed the modular curve known as $X_{11}$.
Since $X$ has genus $1$, the map $\phi_b : X\rightarrow J$ is 
an isomorphism mapping $b=[0:1:0]$ onto the origin. 
So the affine coordinate $x$ and $y$ induce algebraic  functions
$x\circ \phi_b^{-1}$ and $y\circ \phi_b^{-1}$ on $J$.
There are $11^2=121$ points of $11$ torsion in $J$
and $0$ is one of them. We will be interested in the values taken
by  $x\circ \phi^{-1}_b$  at the remaining $120$ points of $11$-torsion.
On can check that $x\circ \phi^{-1}_b$ takes the same value at two opposite
points. So there only remain $60$ values of interest.
These are algebraic numbers and they form a single orbit under
the action of the Galois group $\Gal(\bQ/\QQ)$. So it is natural
to consider their annihilating polynomial
\begin{subeqn}\label{eq:divpol}
H(T)=\prod_{0\not = a\in J[11]/{\pm 1}}\left(T-x(\phi^{-1}(a))\right).
\end{subeqn}
This is an irreducible polynomial in $\QQ[T]$. Computing such
polynomials
is a cornerstone in the algorithmic of modular forms and Galois
representations.

\subsection{An algebraic approach}\label{sec:aaa}

The polynomial in equation~(\ref{eq:divpol}) is known as the
$11$-th division polynomial $\psi_{11}$ of the genus one curve $X$.
For every $k\ge 1$ one can define the $k$-th division polynomial
$\psi_k(T)$ to be the annihilating polynomial of the $x$-coordinates
of all non-zero $k$-torsion points on $X$.  These polynomials can be computed
using  recursion formulae \cite[Section 3.6]{enge} 
\cite[Exercise 3.7]{silverman} that follow from the simple algebraic
form
of the addition law on $X$. Using these recursion
formulae we find
\[
H(T)=T^{60}-20T^{59}+112T^{58}+1855T^{57}+\cdots    + 1321T^4 - 181T^3
+ 22T^2 - 2T + 1/11.
\]
So we have an efficient algebraic method to  compute $H(T)$.
We will explain in section~\ref{sec:ctd2} why it seems difficult to us
to generalize  this algebraic method to curves of higher genus.

\subsection{Using complex approximations}\label{sec:uca}

In this section we compute complex approximations of the coefficients
of $H(T)$. We also explain how one can  deduce the exact value
of these coefficients from a sharp enough complex approximation.
We have seen in sections~\ref{sec:jacinv} and~\ref{sec:psi} how 
to invert the map $\phi_b$. Given a point $a$ in the torus
$\CC^g/\Lambda$ we can compute a complex approximation of some
reduced 
divisor $P_a-gb$ such that $\phi_b (P_a-bg)=a$. Since here the genus
is one, $P_a$ consists of a single point on $X$, and it is uniquely 
defined.  In case $a=(\Omega_1+\Omega_2)/11$ we already found that
the $x$-coordinate
$x(P_a)$ of $P_a$ is 
\[6.796539142094915911068237206 - 7.525908029899464321854796862i\]
up to an error of $10^{-27}$.  We let $a$ run over the $60$ elements
in $(J[11]-\{0\})/\pm 1 $ and compute the $60$ corresponding values of
$x(P_a)$ with the same accuracy. We then compute their sum and find it
is equal to $20$ up to an error of $10^{-25}$. This suggests that the
coefficient of $T^{59}$ in $H(T)$ is $-20$. In order to turn this
heuristic into a proof, we need some information about the
coefficients of $H(T)$. We know that these coefficients are rational
numbers. We need an upper bound on their {\it height}.  The {\it
  height} of a rational number is the maximum of the absolute values
of its numerator and denominator. We explain in the next
section~\ref{sec:recover} how a good approximation and a good bound on
the height suffice to characterise and compute a rational number.  In
case $X$ is the curve given in equation~(\ref{eq:X11}) an upper bound
on the height of the coefficients of $H(T)$ can be proved by
elementary means. For example we know that the denominator of these
coefficients is either $1$ or $11$.  In case $X$ is a modular curve,
similar bounds will be necessary.  These bounds have been proved by
the second author in collaboration with de Jong in \cite{Ed-Jo1} and
\cite{Ed-Jo2}, using Arakelov theory and arithmetic geometry together
with a result of Merkl in \cite{Merkl} on upper bounds for Green
functions.

All the coefficients of $H(T)$ are computed in the same way.
They are symmetric  functions of the  $x(P_a)$, so we
can compute sharp approximations for them. We deduce their
exact values using an a priori  bound on their
height.

\subsection{Recovering a rational number from a good
  approximation}\label{sec:recover}

In the previous section~\ref{sec:uca} we claimed that 
a rational number $x=a/b$ can be recovered from a 
sharp enough complex approximation, provided we have an a priori
bound on the height of $x$.
We recall that the
height of a rational number $a/b$, with $a$
and $b$ integers that are relatively prime, is $\max\{|a|,|b|\}$.
The rational number $x=a/b$ is
known if we know an upper bound $h$ for its height and an
approximation $y$ of it (in $\RR$, say), with
$|x-y|<1/(2h^2)$. Indeed, if $x'=a'/b'$ also has height at most~$h$,
and $x'\neq x$, then
\[
|x-x'| = \left|\frac{a}{b}-\frac{a'}{b'}\right| = 
\left|\frac{ab'-ba'}{bb'}\right| \geq \frac{1}{|bb'|} \geq 1/h^2.
\]
We also note that there are good algorithms to deduce $x$ from such a
pair of an approximation~$y$ and a bound~$h$, for example, by using
continued fractions, as we will now explain. 

In practice we will use rational approximations $y$ of~$x$. Every
rational number $y$ can be written uniquely as
\[
[a_0,a_1,\ldots,a_n] = a_0+\cfrac{1}{a_1+
\cfrac{1}{\ddots\genfrac{}{}{0pt}{0}{}{a_{n-1}+\cfrac{1}{a_n}}}}\;,
\]
where $n\in\ZZ_{\geq0}$, $a_0\in\ZZ$, $a_i\in\ZZ_{>0}$ for all $i>0$,
and $a_n>1$ if $n>0$. To find these $a_i$, one defines $a_0:=\lfloor
y\rfloor$ and puts $n=0$ if $y=a_0$; otherwise, one puts
$y_1:=1/(y-a_0)$ and $a_1=\lfloor y_1\rfloor$ and $n=1$ if $y_1=a_1$,
and so on. The rational numbers $[a_0,a_1,\ldots,a_i]$ with
$0\leq i\leq n$ are called the \emph{convergents} of the continued
fraction of~$y$. Then one has the following well-known result
(see Theorem~184 from~\cite{Hardy-Wright}).
\begin{subproposition}\label{prop_cont_frac}
Let $y$ be in $\QQ$, $a$ and $b$ in~$\ZZ$ with $b\neq 0$, and
\[
\left|\frac{a}{b}-y\right|<\frac{1}{2b^2} .
\]
Then $a/b$ is a convergent of the continued fraction of~$y$. 
\end{subproposition}

\section{Computing torsion points  II}\label{sec:ctd2}

In this section we describe how we compute the fields of definition of
certain torsion points in Jacobians of modular curves. We recommend
\cite{Diamond-Shurman} to those who are interested in an introduction
to the theory of modular forms.

Let $\SL_2(\ZZ)$ denote the group of $2$ by $2$ matrices with
coefficients in $\ZZ$ and with determinant one. It acts on the complex
upper half plane $\HH$ via fractional linear transformations
\begin{eqn}
\left(\begin{matrix}a & b\\ c & d\end{matrix}\right)\cdot z = 
\frac{az+b}{cz+d}.
\end{eqn}
The standard fundamental
domain $F$ for $\SL_2(\ZZ)$ acting on $\HH$ (see Figure~\ref{figbas1})
consists of the $z$ with $|z|\geq1$ and $|\Re(z)|\leq1/2$. 
\begin{figure}
\[
\begin{tikzpicture}[scale=2,thick]
\fill[black!20!white] (0.5,2) -- (0.5,0.866) arc (60:120:1cm) 
-- (-0.5,2)-- cycle;
\draw (-2,0) -- (2,0);
\draw (-0.5,2) -- (-0.5,0) node[anchor=north]{$-1/2$};
\draw (0.5,2) -- (0.5,0) node[anchor=north]{$1/2$};
\draw  (1,0) arc (0:180:1cm) node[anchor=north]{$-1$};
\draw (1,0) -- (1,0) node[anchor=north]{$1$};
\draw (0,1.2) -- (0,1.2) node[anchor=south]{$F$};
\end{tikzpicture}
\]
\caption{Standard fundamental domain $F$ for $\SL_2(\ZZ)$ acting on $\HH$}
\label{figbas1}
\end{figure}
It is not bounded, hence not compact. Viewing $\HH$ as the open
northern hemisphere in $\PP^1(\CC)=\CC\cup\{\infty\}$, with boundary
the equator $\PP^1(\RR)$, we see that the closure $\ol{F}$ of $F$ in
$\PP^1(\CC)$ is the union of $F$ and the point~$\infty$. 

For every prime number $\ell$ we let $\Gamma_\ell$ denote the subset of
$\SL_2(\ZZ)$ consisting of the 
$(\begin{smallmatrix}a & b\\ c & d\end{smallmatrix})$ with $c$, $a-1$
and $d-1$ divisible by~$\ell$. Then $\Gamma_\ell$ is a subgroup of
$\SL_2(\ZZ)$, of index $\ell^2-1$. We assume that $\ell\geq 5$ from
now on. Then the action of $\Gamma_\ell$ on~$\HH$ is free. Each $z$ in
$\HH$ has a neighbourhood $U$
such that all $\gamma U$ for $\gamma$ in $\Gamma_\ell$ are
disjoint. The quotient $\Gamma_\ell\backslash\HH$ is therefore a
Riemann surface that we denote by~$Y_\ell$, and the quotient map
$\HH\to Y_\ell$ is a covering map, that is, each point $y$ in $Y_\ell$
has an open neighbourhood $U$ such that the inverse image of $U$ in
$\HH$ is the disjoint union of copies of $U$, indexed by the inverse
image of~$y$.

The Riemann surface $Y_\ell$ is not compact. A fundamental domain
$F_\ell$ in $\HH$ for $\Gamma_\ell$ can be gotten as the union of the
$\gamma F$, where $\gamma$ ranges over a set of representatives of
$\Gamma_\ell\backslash\SL_2(\ZZ)/\{1,-1\}$. Such a set consists of
$(\ell^2-1)/2$ elements and it can easily be found.  We can compactify
$Y_\ell$ to a compact Riemann surface $X_\ell$ by adding $\ell-1$
points, called \emph{cusps}, the points of $\PP^1(\RR)$ that lie in
the closure of $F_\ell$ in $\PP^1(\CC)$. These points lie in fact in
$\PP^1(\QQ)$ and can easily be written down. All this leads to an
explicit topological and analytic description of~$X_\ell$. It is
covered by coordinate disks around the cusps. For example, the
function
\begin{eqn}
q\colon\HH\to\CC, \quad z\mapsto e^{2\pi iz},
\end{eqn}
restricted to the set of $z$ with $\Im(z)>1/\ell$, induces a
coordinate on a disk in $X_\ell$ around the cusp~$\infty$. Indeed, the
image under $q$ of this region is the punctured disk of radius
$e^{-2\pi/\ell}$ around~$0$, and the cusp $\infty$ fills the
puncture. The genus $g_\ell$ of $X_\ell$ is equal to
$(\ell-5)(\ell-7)/24$. For $\ell=11$ the genus is~$1$, and indeed,
$X_{11}$ is the elliptic curve $X_{11}$ given by
equation~(\ref{eq:X11}).

It is of course a miracle that such an analytically defined Riemann
surface as $X_{11}$ is defined over $\QQ$, that is, can be
described as a curve in a projective space given by a equations with
coefficients in~$\QQ$. But this is true for all~$\ell$, and it is
explained as follows, for $\ell>13$. The theory of modular forms
gives that the $\CC$-vector spaces $\Omega^1(X_\ell)$ of holomorphic
differentials on $X_\ell$ have bases consisting of $1$-forms $\omega$
whose pullback to $\HH$ is of the form
$(\sum_{n\geq1}a_nq^n){\cdot}(dq)/q$ with all $a_n$ in~$\ZZ$. Quotients of
such $\omega$ and $\omega'$ in $\Omega^1(X_\ell)$ then provide
sufficiently many rational functions on $X_\ell$ to embed it into a
projective space, such that the image is given by homogeneous
polynomial equations with coefficients in~$\QQ$.

We let $J_\ell$ denote the Jacobian variety of $X_\ell$. It is also
defined over~$\QQ$, as well as its group law. This means that the
group law is described by quotients of polynomials with coefficients
in~$\QQ$. Therefore, for all $P$ and $Q$ in $J_\ell$ and for each
$\sigma$ in $\Aut(\CC)$, the automorphism group of the field~$\CC$, we
have $\sigma(P+Q)=\sigma(P)+\sigma(Q)$. For each integer $m\geq1$ the
kernel $J_\ell[m]$ of the multiplication by $m$ map is finite (it
consists of $m^{2g_\ell}$ elements) and preserved by the action of
$\Aut(\CC)$. This implies that all $P$ in $J_\ell[m]$ have coordinates
in the algebraic closure $\bQ$ of $\QQ$ in $\CC$, that is, for each
rational function $f$ on $J_\ell$ that is defined over $\QQ$ and has
no pole at~$P$, the value $f(P)$ of $f$ at $P$ is in~$\bQ$. The
analytic description above of $X_\ell$ gives us an analytic
description of~$J_\ell$.

We are interested in certain subgroups $V_\ell$ of the $\ell$-torsion
subgroup $J_\ell[\ell]$ of $J_\ell$ that are invariant under the
Galois group $\Gal(\bQ/\QQ)$ and consist of $\ell^2$ elements. These
$V_\ell$ can be described explicitly and efficiently in terms of
certain operators called Hecke operators on the first homology group
of~$X_\ell$. The whole point is to understand them algebraically, with
their $\Gal(\bQ/\QQ)$-action.

The subgroup $V_\ell$ defines a commutative $\QQ$-algebra $A_\ell$ of
dimension $\ell^2$ as $\QQ$-vector space, the coordinate ring of
$V_\ell$ over~$\QQ$. This algebra $A_\ell$ consists of the functions
$f\colon V_\ell\to\bQ$ with the property that for all $\sigma$ in
$\Gal(\bQ/\QQ)$ and all $P$ in $V_\ell$ we have $f(\sigma
(P))=\sigma(f(P))$. Addition and multiplication are pointwise. Each
$f_\ell$ in $A_\ell$ with the property that the $f_\ell(P)$ are all
distinct is a generator, and $A_\ell$ is then given as
$\QQ[T]/(H_{f_\ell})$, with
\[
H_{f_\ell} = \prod_{P\in V_\ell}(T-f_\ell(P)) \quad\text{in $\QQ[T]$.}
\]

A direct approach for computing $A_\ell$ or $H_{f_\ell}$
algebraically, as in Section~\ref{sec:aaa} in the case of the division
polynomial $\psi_{11}$, is very unlikely to succeed in time polynomial
in~$\ell$, because in the case of $V_\ell$ one has to work with the
algebraic variety $J_\ell$, whose dimension grows quadratically
with~$\ell$. Writing down polynomial equations with coefficients in
$\QQ$ for $J_\ell$ and $V_\ell$ is probably still possible, in time
polynomial in~$\ell$. But computing a $\QQ$-basis of $A_\ell$ from the
equations in a standard way uses Groebner basis methods, which, as far
as we know, take time exponential or even worse in the number of
variables, that is, exponential or worse in~$\ell$.

For this reason we replace, in \cite{book}, exact computations by
approximations. There are then two problems to be dealt with. The
first is to show that $f_\ell$ can be chosen so that the logarithm of
the height of the coefficients of $H_{f_\ell}$, that is, the number of
digits of their numerator and denominator, does not grow faster than a
power of~$\ell$. This problem is solved in \cite{Ed-Jo1}, \cite{Merkl}
and \cite{Ed-Jo2}, using arithmetic algebraic geometry and analysis on
Riemann surfaces. The second problem is to show that for the same
choice of $f_\ell$, the values $f_\ell(P)$ at all $P$ in $V_\ell$ can
be approximated in $\CC$ with a precision of $n$ digits in time
polynomial in $n+\ell$. This is done in~\cite{Couveignes2}. The
chapters \cite{Bosman1} and \cite{Bosman2} contain real computations
using the method of Section~\ref{sec:jacinv}, for prime numbers
$\ell\leq 23$.

Let us now explain how we choose $f_\ell$ (up to some technicalities;
the precise setup is given in~\cite[\S8.2]{Ed3}) and say some words
about the approximation of the~$f_\ell(P)$. Standard functions on
Jacobian varieties such as $J_\ell$ are theta functions. But a problem
is that these are usually given as power series in $g_\ell$ variables,
and as $g_\ell$ grows this can make the number of terms that must be
evaluated for a sufficiently good approximation grow exponentially
in~$\ell$. In other words, we know no method to approximate their
values fast enough (of course, it is not excluded that such methods do
exist). Our solution is to transfer the problem from $J_\ell$ to
$X_\ell^{g_\ell}$, via the Abel-Jacobi map. We choose $h_\ell$ a
suitable non-constant rational function on $X_\ell$, defined over
$\QQ$, of small degree and with small coefficients. Then we take as
origin a suitable divisor of degree $g_\ell$ on $X_\ell$, defined
over~$\QQ$. This divisor is carefully chosen in~\cite{Ed3} to have the
following property: for each $P$ in $V_\ell$ there is a {\it unique}
effective divisor $Q_P = Q_{P,1}+\cdots+Q_{P,g_\ell}$ on $X_\ell$,
such that its image under the Abel-Jacobi map is~$P$. Then we define
$f_\ell(P)=h_\ell(Q_{P,1})+\cdots+h_\ell(Q_{P,g_\ell})$. Rather
magically, the problem of power series in many variables has
disappeared. The function $h_\ell$ is locally given by a power series
in {\it one} variable. We evaluate it at each $Q_{P,i}$ separately.
The Abel-Jacobi map (see Section~\ref{sec:afcaj}) is given by a sum of
$g_\ell$ integrals of $g_\ell$-tuples of holomorphic $1$-forms in one
variable. The analytic description above of $X_\ell$ and $J_\ell$
should make it clear that the Abel-Jacobi map and the function
$h_\ell$ can be well approximated with standard tools. That means that
the only remaining problem is the inversion of the Abel-Jacobi map,
that is, the approximation of the divisors $Q_P$, but that was
discussed and solved in Sections~\ref{sec:afcaj}
and~\ref{sec:psi}. The main result obtained in \cite{book} is the
following theorem.
\begin{theorem}\label{thmcompVl}
There is a deterministic algorithm that on input a prime number
$\ell\geq 11$ computes the $\QQ$-algebra $A_\ell$ in time polynomial
in~$\ell$.
\end{theorem}

A {\it probabilistic} algorithm for computing $A_\ell$ is also given
in~\cite{book}. It relies on $p$-adic approximations rather than
complex approximations. In~\cite{Couveignes3} it is explained how such
$p$-adic approximations can be computed efficiently.  From a
theoretical point of view, a probabilistic algorithm is not quite as
satisfactory as a deterministic one. From a practical point of view,
it is just as good. In our case the probabilistic algorithm has a
simpler proof than the deterministic one. And Peter Bruin
\cite{Bruin1,Bruin2} has been able to generalize it to a much wider
class of $V_\ell$-like modular spaces. Finding a similar
generalization for the deterministic algorithm is an open problem at
this time.

\section{Applications and  open questions}\label{sec:avoq}
The main motivation for all the work done in \cite{book} is the
application in number theory to the fast computation of coefficients
of modular forms. Instead of attempting to present this in the
most general case we give two examples: Ramanujan's $\tau$-function,
and powers of Jacobi's $\theta$-function. 

Recall that $q\colon \HH\to\CC$ is the function $z\mapsto e^{2\pi iz}$. 
The \emph{discriminant modular form}~$\Delta$ is the holomorphic
function on $\HH$ given by the converging infinite product
\begin{eqn}\label{eqn5.1}
\Delta = q\prod_{n\geq1}(1-q^n)^{24}.
\end{eqn}
The holomorphic function $\Delta$ has a power series expansion in $q$, 
\begin{eqn}
\Delta = \sum_{n\geq1}\tau(n)q^n, 
\end{eqn}
whose coefficients, which are integers, define \emph{Ramanujan's
  $\tau$-function}. It can be shown that for all
$(\begin{smallmatrix}a & b\\ c & d\end{smallmatrix})$ in $\SL_2(\ZZ)$,
  and for all $z$ in $\HH$, we have
\begin{eqn}\label{eqn5.3}
\Delta\left(\frac{az+b}{cz+d}\right) = (cz+d)^{12}\Delta(z).    
\end{eqn}
Functions $f\colon\HH\to\CC$ that are given by a power series
$\sum_{n\geq 1}a_n(f)q^n$ with this symmetry under the action of
$\SL_2(\ZZ)$ on $\HH$ with the exponent $12$ replaced by an integer
$k$ are called cuspidal modular forms of weight $k$ on
$\SL_2(\ZZ)$. The complex vector spaces $S(\SL_2(\ZZ),k)$ of cuspidal
modular forms of weight $k$ are finite dimensional. The dimension
grows roughly as~$k/12$. More precisely, for $k<12$ the space
$S(\SL_2(\ZZ),k)$ is zero, and $S(\SL_2(\ZZ),12)$ is one-dimensional,
generated by~$\Delta$. The fact that each $g$ in $\GL_2(\QQ)$ with
$\det(g)>0$ acts on $\HH$ and normalises $\SL_2(\ZZ)$ up to finite
index leads to operators $T_{k,g}$ on the $S(\SL_2(\ZZ),k)$. These
operators are named after Hecke. For each integer $n\geq1$ there is an
operator $T_{k,n}$; for $n$ prime it is the one induced by the matrix
$(\begin{smallmatrix}n&0\\0&1
\end{smallmatrix})$ 
and for general $n$ it is a bit more complicated. As the space
$S(\SL_2(\ZZ),12)$ is one-dimensional, each $T_{12,n}$ acts on it as
multiplication by a scalar. This scalar turns out to be the
coefficient $\tau(n)$ of $q^n$ in the power series of~$\Delta$. Well
known relations between the Hecke operators imply relations between
the $\tau(n)$ that are summarised in the identity of Dirichlet series,
for $s$ in $\CC$ with real part $\Re(s)$ large enough:
\begin{eqn}
\sum_{n\geq1} \tau(n)n^{-s} =  \prod_p(1-\tau(p)p^{-s}+p^{11}p^{-2s})^{-1}.
\end{eqn}
Here the product is over all prime numbers, and both sides converge
for $\Re(s)>13/2$. In fact, it is a famous theorem of Deligne
(\cite{Deligne1} and \cite{Deligne2}) that for all primes $p$ one has
\begin{eqn}\label{eq:ineq}
|\tau(p)|\leq 2p^{11/2},
\end{eqn} as conjectured by Ramanujan. 

The identity of Dirichlet series shows that for $n$ and $m$ with
greatest common divisor $1$ we have $\tau(nm)=\tau(n)\tau(m)$, and
that for $p$ prime and $r\geq2$ we have
$\tau(p^r)=\tau(p)\tau(p^{r-1})-p^{11}\tau(p^{r-2})$. Therefore, the
computation of $\tau(n)$ is reduced to that of the $\tau(p)$ for $p$
dividing~$n$. We can now state one of the main theorems of \cite{book}.

\begin{theorem}\label{theorem:tau}
There is a deterministic algorithm that on input an integer $n\geq1$
together with its factorisation into prime factors computes $\tau(n)$
in time polynomial in $\log n$.
\end{theorem}
Before this result, the fastest known algorithms to compute $\tau(n)$
took time exponential in $\log n$. For example, if one computes the
product in (\ref{eqn5.1}) up to order $n$ by multiplying the necessary
factors, then one spends time at least linear in~$n$. To prove the
theorem, it suffices to show that for $p$ prime $\tau(p)$ can be
computed in time polynomial in $\log p$. This will be done using
Theorem~\ref{thmcompVl}.

The fact that modular forms have an enormous amount of symmetry as in
(\ref{eqn5.3}) is certainly powerful, but it does not suffice at this
point. What is needed is Galois symmetry, which is also behind
Deligne's famous result mentioned above. A lot could be said on this,
but this is not an appropriate place for that. 

In a nutshell: modular forms give elements in de Rham cohomology of
complex algebraic varieties defined over $\QQ$, and the singular
homology with torsion coefficients $\ZZ/\ell\ZZ$ of those complex
varieties can be defined algebraically (Grothendieck, Artin, Deligne)
and therefore has an action by~$\Gal(\bQ/\QQ)$. 

For example, $\Delta$ gives rise, for every prime integer
$\ell\ge 11$,  to a certain subgroup~$V_\ell$ of the $\ell$-torsion
of the Jacobian $J_\ell$  of $X_\ell$. This subgroup $V_\ell$
has cardinality $\ell^2$.
 For $\ell\neq p$
the image of $\tau(p)$ in $\ZZ/\ell\ZZ$ is determined by the action of
$\Gal(\bQ/\QQ)$ on this~$V_\ell$. The addition map $V_\ell\times
V_\ell\to V_\ell$ induces a $\QQ$-algebra morphism called co-addition
from $A_\ell$ to $A_\ell\otimes A_\ell$, that is, from
$\QQ[T]/(H_{f_\ell})$ to
$\QQ[T_1,T_2]/(H_{f_\ell}(T_1),H_{f_\ell}(T_2))$. Computing $\tau(p)$
modulo $\ell$ (for $p\neq\ell$) is then done by reducing $A_\ell$ with
its co-addition modulo $p$ and computing on this reduction
$A_{\ell,p}$ a certain relation between the co-addition and the
Frobenius map that sends $a$ in $A_{\ell,p}$ to~$a^p$, just as in
Schoof's algorithm for elliptic curves (see Section~1.2
of~\cite[\S1.2]{Ed1}). For more detail the interested reader is referred to
Section~2.4 of \cite[\S2.4]{Ed2} and the references therein. The point is
that this advanced machinery can actually be used for computing
$\tau(p)\bmod \ell$ in time polynomial in $\log p$ and $\ell$. 

In order to recover the actual value of $\tau(p)$ as an integer, we
compute $\tau(p)$ modulo several small primes $\ell$. If the product
of these small primes is bigger than $4p^{5.5}$ then 
we deduce $\tau(p)$ using inequality~(\ref{eq:ineq}) and Chinese
remainder theorem \cite[1.3.3]{Cohen}.  

We now come to our second example: the classical question in how many
ways a positive integer $n$ can be written as a sum of $d\geq 1$
squares of integers. Let us write $r_d(n)$ for this number, that is,
$r_d(n)=\#\{x\in\ZZ^d : x_1^2+\cdots+x_d^2=n\}$. Then $r_d(n)$ is the
coefficient of $q^n$ in the formal power series~$\theta_d$, with
\begin{eqn}
\theta_d = \sum_{n\geq 0}r_d(n)q^n = \sum_{x\in\ZZ^d}q^{x_1^2+\cdots+x_d^2}=
\left(\sum_{x_1\in\ZZ}q^{x_1^2}\right)\cdots 
\left(\sum_{x_d\in\ZZ}q^{x_d^2}\right) =
\theta_1^d \quad\text{in $\ZZ[[q]]$}.
\end{eqn}
The formal power series $\theta_1$ defines a holomorphic complex
function on the complex upper half plane $\theta\colon\HH\to\CC$ by
viewing $q$ as the function $q\colon z\mapsto e^{2\pi iz}$. Poisson's
summation formula then shows that for all $z\in\HH$ we have
\begin{eqn}
\theta(-1/4z) = (-2iz)^{1/2}\theta(z),
\end{eqn}
where the square root is continuous and positive for $z$
in~$i{\cdot}\RR_{>0}$. This functional equation for $\theta$, together
with the obvious one $\theta(z+1)=\theta(z)$, imply that $\theta$ is a
modular form of weight $1/2$, and therefore that $\theta_d$
(interpreted as a function on~$\HH$) is a modular form of
weight~$d/2$.

This fact is the origin of many results concerning the
numbers~$r_d(n)$. The famous explicit formulas for the $r_d(n)$ for
even $d$ up to $10$ due to Jacobi, Eisenstein and Liouville (see
~\cite{Milne2} and Chapter~20 of~\cite{Hardy-Wright}) owe their
existence to it. In order to state these formulas, let $\sum_{d|m}$
denote summation over the positive divisors $d$ of~$m$, with the
convention that there are no such $d$ if $m$ is not an integer, and
let $\chi\colon\ZZ\to\CC$ be the map that sends $n$ to $0$ if $n$ is
even, to $1$ if $n$ is of the form $4m+1$ and to $-1$ if $n$ is of the
form $4m-1$. Then we have:
\[
\begin{aligned}
r_{2}(n) & = 4\sum_{d|n}\chi(d), \\
r_{4}(n) & =8\sum_{2\nmid d|n}d + 16\sum_{2\nmid d|(n/2)}d, \\
r_{6}(n) & = 16\sum_{d|n}\chi\left(\frac{n}{d}\right)d^2 
 -4\sum_{d|n}\chi(d)d^2, \\
r_{8}(n) & = 16\sum_{d|n}d^3 - 32\sum_{d|(n/2)}d^3 +
256\sum_{d|(n/4)}d^3, \\
r_{{10}}(n) & = \frac{4}{5}\sum_{d\mid n}\chi(d)d^4 + 
\frac{64}{5}\sum_{d\mid n}\chi\left(\frac{n}{d}\right)d^4  
+\frac{8}{5}\sum_{d\in\ZZ[i],\, |d|^2=n}d^4.
\end{aligned}
\]
In the last formula, $\ZZ[i]$ is the set of Gaussian integers $a+bi$
in $\CC$ with $a$ and $b$ in~$\ZZ$.  

Using these formulas, the numbers $r_d(n)$ for $d$ in $\{2,4,6,8,10\}$
can be computed in time polynomial in $\log n$, if $n$ is given with
its factorisation in prime numbers. This is not the case for the
formulas that were found a bit later by Glaisher for $r_d(n)$ for some
even $d\geq12$. We give the formula that he found for $d=12$, as
interpreted by Ramanujan:
\begin{eqn}
r_{12}(n) = 8\sum_{d|n}d^5 -512\sum_{d|(n/4)}d^5 + 16a_n,
\quad\text{where}
\quad \sum_{n\geq 1}a_nq^n = q\prod_{m\geq 1}(1-q^{2m})^{12}.
\end{eqn}
Computing $a_n$ by multiplying out the factors $1-q^{2m}$ up to order
$n$ takes time at least linear in $n$, hence exponential in $\log
n$. We know of no direct way to compute the $a_n$ in time polynomial
in $\log n$, even if $n$ is given with its factorisation. However,
$\sum_{n\geq1}a_nq^n$ is a modular form, and therefore we \emph{can}
compute $a_n$ in time polynomial in $\log n$, if $n$ is given with its
factorisation, \emph{via} the computation of Galois representations.
The same is true for the $r_d(n)$ for all even~$d$. The explicit
formulas for $d\leq 10$ correspond precisely to the cases where the
Galois representations that occur are of dimension one, whereas for
$d\geq 12$ genuine two-dimensional Galois representations always
occur, as proved by Ila Varma in her master's thesis~\cite{Varma}. 

We conclude that from an algorithmic perspective the classical problem
of computing the $r_d(n)$ for even $d$ and $n$ given with its
factorisation into primes is solved for \emph{all} even~$d$. The
question as to the existence of \emph{formulas} has a negative answer,
but for \emph{computations} this does not matter.

\paragraph{Open questions}
Finally, we should point out that
the algorithms in theorems~\ref{thmcompVl} and~\ref{theorem:tau},
despite their polynomial time complexity,
are  not so practical  at present. 
However, Bosman's computation of the $V_\ell$
associated with  $\Delta$ for  $\ell$ in $\{13, 17, 19\}$ enabled him
to further study
Lehmer's conjecture on the values of $\tau(n)$ modulo $n$.
See Lygeros and Rozier~\cite{Ly-Ro} for a more
classical experimental approach. A challenge for the near future is to
design and implement a practical variant of these algorithms.

\end{document}